\documentclass[reqno,10pt]{amsart}
\usepackage{amscd,amssymb,verbatim,array}
\usepackage{hyperref}

\numberwithin{equation}{section} \theoremstyle{plain}

\newtheorem{theorem}{Theorem}[section]

\newtheorem{proposition}[theorem]{Proposition}

\newtheorem{definition}[theorem]{Definition}
\newtheorem{conjecture}[theorem]{Conjecture}

\theoremstyle{definition}

\theoremstyle{remark}

\numberwithin{equation}{section}

\newcommand{\B}{\operatorname{B}}
\newcommand{\Eul}{\operatorname{Eul}}
\newcommand{\even}{\operatorname{even}}

\newcommand{\Mod}{\operatorname{mod}}
\newcommand{\rk}{\operatorname{rank}}

\begin{document}

\title{REFINED ANALYTIC TORSION: COMPARISON THEOREMS AND EXAMPLES}

\author{RUNG-TZUNG HUANG}

\email{huang.r@neu.edu}

\address{Department of Mathematics, Northeastern University, Boston, MA  02115, USA}



\begin{abstract}
Braverman and Kappeler introduced a refinement of the Ray-Singer analytic
torsion associated to a flat vector bundle over a closed
odd-dimensional manifold. We study this notion and improve the
Braverman-Kappeler theorem comparing the refined analytic torsion
with Farber-Turaev refinement of the combinatorial torsion. Using
this result we establish, modulo sign, the Burghelea-Haller
conjecture, comparing their complex analytic torsion with
Farber-Turaev torsion in the case, when the flat connection can be
deformed in the space of flat connections to a Hermitian connection.
We then compute the refined analytic torsion of lens spaces and
answer some of the questions posed in \cite[Remark 14.6.2]{BK1}.

\end{abstract}

\maketitle


\section{Introduction}
Let $M$ be a closed oriented odd dimensional manifold and let $E$ be
a complex vector bundle over $M$ endowed with a flat connection
$\nabla$. In a series of papers, \cite{BK1,BK2,BK3,BK4}, M.
Braverman and T. Kappeler defined and studied a nonzero element
\[
\rho_{\operatorname{an}}(\nabla) \in
\operatorname{Det}\big(H^{\bullet}(M,E)\big)
\]
of the complex determinant line
$\operatorname{Det}\big(H^{\bullet}(M,E)\big)$ of the cohomology
$H^{\bullet}(M,E)$ of $M$ with coefficients in the complex vector
bundle $E$. They called this element {\em refined analytic torsion}.
It can be viewed as an analytic analogue of the refinement of the
Reidemeister torsion due to Turaev \cite{T1,T2} and, more generally,
to Farber and Turaev \cite{FT1,FT2}. Recall that the Farber-Turaev
torsion $\rho_{\varepsilon,o}(\nabla)$ depends on the Euler
structure $\varepsilon$, the cohomology orientation $o$, and the
connection $\nabla$.

The following extension of the Cheeger-M\"{u}ller theorem \cite{C,M}
was proven in \cite[Theorem 5.11]{BK4}: For each connected component
$C$ of the space $\text{Flat}(E)$ of flat connections on $E$, there
exists a constant $\theta^C \in \mathbb{R}$, such that
\begin{equation}\label{E:ratiointrod}
\frac{\rho_{\operatorname{an}}(\nabla)}{\rho_{\varepsilon,o}(\nabla)}
\ =\  e^{i\theta^C} \cdot f_{\varepsilon,o}(\nabla),
\end{equation}
where $f_{\varepsilon,o}(\nabla)$ is a holomorphic function of
$\nabla \in \text{Flat}(E)$, given by an explicit local expression.
Equality \eqref{E:ratiointrod} does not give us any information
about the constant $\theta^C$ and its dependence on $C$,
$\varepsilon$, and $o$. In particular, Braverman and Kappeler posed
the following two questions in \cite[Remark 14.6.2]{BK1},
\subsection*{Question 1} {\em Does the constant\/ $\theta^C$
depend on the connected component\/ $C$ of the space
$\operatorname{Flat}(E)$ of flat connections on $E$?}
\subsection*{Question 2} {\em For which connections\/ $\nabla$ one can find an Euler
structure $\varepsilon$ and the cohomological orientation $o$ such
that\/ $\rho_{\operatorname{an}}(\nabla) =
\rho_{\varepsilon,o}(\nabla)$?}
\medskip

In Section 3 we compute the constant $\theta^C$ for any connected
component $C$ of $\operatorname{Flat}(E)$ which contains a Hermitian
connection. In Section 5 and Section 6 of this paper we compute the
refined analytic torsion of lens spaces and study its relationship
with the cohomological Turaev torsion of lens spaces. Our explicit
calculation for the three-dimensional lens space $L(5;1,1)$ shows
that, in general, {\em the constant\/ $\theta^C$ does depend on the
connected component $C$ of $\operatorname{Flat}(E)$}. This result
provides a positive answer to Question~1 above. (Note that, in the
case of lens spaces, $\operatorname{Flat}(E)$ is discrete and
coincides with the space of acyclic Hermitian connections. Hence,
connected components $C$ of $\operatorname{Flat}(E)$ are one-element
subsets).

We then compute the quotient of the refined analytic torsion and
cohomological Turaev torsion of the five-dimensional lens space
$L(3; 1,1,1)$. In this case we show that {\em for all connections,
all Euler structures and all cohomological orientations, the
cohomological Turaev torsion and the refined analytic torsion are
not equal}. This provides a partial answer to Question~2 above.

In \cite{BH1, BH2} Burghelea and Haller defined a complex valued
quadratic form, referred to as {\em complex Ray-Singer torsion}. This torsion is defined for a complex flat vector bundle over a closed
manifold of arbitrary dimension, provided that the complex vector
bundle admits a non-degenerate complex valued symmetric bilinear
form $b$. Burghelea and Haller, \cite[Conjecture 5.1]{BH2}$^\mathrm{1}$, see also
Conjecture~\ref{C:conj} below, conjectured that the complex Ray-Singer
torsion is roughly speaking equal to the square of the Farber-Turaev
torsion and established the conjecture in some non-trivial
situations. \footnotetext[1]{The author was informed by Burghelea and Haller that they can prove \cite[Conjecture 5.1]{BH2} for all odd dimensional manifolds up to sign.} Braverman and Kappeler, \cite{BK5}, expressed the
Burghelea-Haller complex Ray-Singer torsion in terms of the square
of the refined analytic torsion $\rho_{\operatorname{an}}(\nabla)$
and the eta invariant $\eta(\nabla)$. In particular, they proved a
weak version of the Burghelea-Haller conjecture$^\mathrm{2}$.\footnotetext[2]{In \cite[Theorem 5.9]{BH2}, Burghelea and Haller also obtain the same result and generalize it to the even dimensional case using the idea of Braverman and Kappeler in \cite{BK5}.} In Section 4 we improve this result for the case when $\nabla$ belongs to a
connected component of the space of flat connections on the
associated complex vector bundle $E$ which contains a Hermitian
connection. Our result establishes, modulo sign, the Burghelea-Haller conjecture
for this case.

This paper is organized as follows. In Section 2 we recall the
definitions and properties of refined analytic torsion from
\cite{BK1,BK2,BK3}. In Section 3 we studied the comparison theorem
of the refined analytic torsion and the cohomological Farber-Turaev
torsion from \cite[Theorem 5.11]{BK4} and present the formula of the
constant $\theta^C$. In Section 4 we present our result about the
Burghelea-Haller conjecture. In Section 5 we compute the refined
analytic torsion of lens spaces. In Section 6 we compute the Turaev
torsion of lens spaces. In the end of Section 6 we calculate the
constant $\theta^C$ in the case of the three-dimensional lens space
$L(5; 1,1)$ and the quotient of the refined analytic torsion and
cohomological Turaev torsion of the five-dimensional lens space
$L(3; 1,1,1)$ and explain how our computation gives answers to
Question~1 and Question~2 above.

\subsection*{Acknowledgment}
The author would like to thank Maxim Braverman for suggesting these
problems as well as for help and encouragement throughout. The author is also
very grateful to Thomas Kappeler and to the referee for very
valuable comments and suggestions.


\section{Refined Analytic Torsion}
Throughout this paper we will assume that $M$ is a closed oriented
manifold of odd dimension $d = 2n-1$ and $E$ is a complex vector
bundle over $M$ endowed with a flat connection $\nabla$. Fix a
Riemannian metric $g^M$ on $M$. In \cite{BK3}, Braverman and
Kappeler defined a non-zero element, cf. \cite[Section 7]{BK3},
\begin{equation}
\rho(\nabla,g^M) \in \operatorname{Det}\big(H^{\bullet}(M,E)\big).
\end{equation}
In general, $\rho(\nabla,g^M)$ might depend on the Riemannian metric $g^M$. Hence they introduced {\em the refined analytic torsion} $\rho_{\operatorname{an}}(\nabla)$, cf. Definition \ref{D:def}, which is a slight modification of $\rho(\nabla,g^M)$ and is independent of $g^M$.

\subsection{The odd signature operator}
The refined analytic torsion is defined in terms of the odd
signature operator which was introduced by Atiyah, Patodi, and
Singer, \cite[p. 44]{APS1}, \cite[p. 405]{APS2}, and, in the more
general setting, by Gilkey, \cite[p. 64-65]{G}. Hence, let us begin
by recalling the definition of this operator.

Let $\Omega^{\bullet}(M,E)$ denote the space of smooth differential forms on $M$ with values in $E$. Fix a Riemannian metric $g^M$ on $M$ and let
$*:\Omega^{\bullet}(M,E)\to \Omega^{d-\bullet}(M,E)$ denote the Hodge $*$-operator. Define the
{\em chirality operator} $\Gamma =  \Gamma(g^M):\Omega^{\bullet}(M,E)\to \Omega^{\bullet}(M,E)$ by the formula
\begin{equation}
\Gamma \omega \,:=\, i^n(-1)^{\frac{k(k+1)}{2}}*\omega, \qquad \omega\in\Omega^k(M,E),
\end{equation}
where $n$ is given as above by $n=\frac{d+1}{2}$. Note that $\Gamma^2=1$.

\begin{definition}
The (even part of the) {\em odd signature operator} $\B_{\even}=\B_{\even}(\nabla,g^M)$ acting on an even form $\omega \in \Omega^{2p}(M,E)$ is defined by the formula
\[
    \B_{\even}\,\omega \ := \ (\Gamma \nabla + \nabla \Gamma) \omega \ \in \ \Omega^{d-2p-1}(\,M,E\,)\,\oplus\,\Omega^{d-2p+1}(\,M,E\,).
\]
\end{definition}

The operator $\B_{\even}$ is an elliptic differential operator, whose leading symbol is symmetric with respect to any Hermitian metric $h^E$ on $E$.


\subsection{The $\eta$-invariant}
Let $\theta$ be an Agmon angle for $\B_{\even}$, see \cite[Definition 3.4]{BK1} or \cite[Definition 6.3]{BK3} for the choice
of this angle. The $\eta$-function of $\B_{\even}$ is defined by the
formula
\[
    \eta_{\theta}(s,\B_{\even}) \ = \
    \sum_{\operatorname{Re}\lambda_k>0}m_k(\lambda_k)^{-s} -
    \sum_{\operatorname{Re}\lambda_k<0}m_k(-\lambda_k)^{-s},
\]
here $\lambda_k$ is the eigenvalue of $\B_{\even}$ and $m_k$ is the algebraic multiplicity of $\lambda_k$. It is known, \cite{G}, that $\eta_{\theta}(s,\B_{\even})$ has a
meromorphic extension to the whole complex plane $\mathbb{C}$ with isolated simple poles, and that it is regular at $0$.

Let $m_+$ (respectively, $m_-$) denote the number of eigenvalues (counted with their algebraic multiplicities) of $\B_{\even}$ on the positive (respectively, negative) part of the imaginary axis. Let $m_0$ denote the algebraic multiplicity of $0$ as an eigenvalue of $\B_{\even}$.

\begin{definition}
The $\eta$-invariant $\eta(\nabla)$ of $\B_{\even}$ is defined by
the formula
\[
  \eta(\nabla)  \ = \ \frac{\eta_{\theta}(\,0,\B_{\even}\,) + m_+ - m_-+m_0}{2}.
\]
\end{definition}
Note that  $\eta(\nabla)$ is independent of the angle $\theta$, cf.
\cite[Subsection 3.10]{BK1}.

\subsection{The refined analytic torsion}
Let $\B_{\text{trivial}} \, = \, \Gamma d+d \Gamma \ : \  \Omega^{\bullet}(M)\to \Omega^{\bullet}(M)$.  Define
\[
  \eta_{\text{trivial}} \ = \  \eta_{\text{trivial}}(g^M) \ = \ \frac{\eta_{\theta}(\,0,\B_{\operatorname{trivial}}\,)}{2}
\]
to be the $\eta$-invariant corresponding to the trivial line bundle $M \times \mathbb{C} \to M$ over $M$.

Recall that the element $\rho(\nabla,g^M) \in \operatorname{Det}\big(H^{\bullet}(M,E)\big)$ of the determinant line of the cohomology of $M$ with coefficients in $E$ was defined in~\cite[Section 7]{BK3}. If the bundle $E$ is acyclic, $\operatorname{Det}\big(H^{\bullet}(M,E)\big)$ is canonically isomorphic to $\mathbb{C}$. In this case, $\rho(\nabla,g^M)$ can be viewed as a complex number, which is equal to the {\em graded determinant} of the operator $\B_{\even}$, cf. \cite[Section 6]{BK1}.

\begin{definition}\label{D:def}
Let $(E,\nabla)$ be a flat vector bundle on $M$. The {\em refined analytic torsion} is the element
\begin{equation}\label{E:def}
   \rho_{\operatorname{an}}(\nabla) \ := \
       \rho(\nabla,g^M) \cdot \exp\Big(\,i\pi \cdot \rk E \cdot \eta_{\operatorname{trivial}}(g^M) \,\Big) \ \in\  \operatorname{Det}\big(H^{\bullet}(M,E),
\end{equation}
where $g^M$ is any Riemannian metric on $M$.
\end{definition}

It is shown in \cite[Theorem 9.6]{BK3} that
$\rho_{\operatorname{an}}(\nabla)$ is independent of the choice of
the metric $g^M$. Note that  when $\dim{M} \equiv 1(\Mod 4)$,
$\eta_{\operatorname{trivial}} = 0$, and, hence, $\rho_{\operatorname{an}}(\nabla) \,=\, \rho(\nabla,g^M)$.



\section{Comparison between the refined analytic torsion and the Farber-Turaev torsion}
We now recall the definition of the canonical involution on the
complex determinant line
$\operatorname{Det}\big(H^{\bullet}(M,E)\big)$ of the cohomology
$H^{\bullet}(M, E)$ of $M$ with coefficients in $E$. We then derive
the formula of the phase of the refined analytic torsion and recall
the formula of the phase of the Farber-Turaev torsion. Then we
compute the constant $\theta^C$ and improve the Braverman-Kappeler
theorem comparing the refined analytic torsion and the Farber-Turaev
torsion.


\subsection{Involution on the determinant line.}
In this subsection we recall the definition of the canonical
involution on the complex line
$\operatorname{Det}\big(H^{\bullet}(M,E)\big)$ from \cite[Subsection
10.1]{BK3}.

Let $M$ be a closed oriented manifold of odd dimension $d=2n-1$ and let $E$ be a flat complex vector bundle over $M$ admitting a flat Hermitian metric $h^E$ and endowed with a flat connection $\nabla$. Then
\[
dh^E(u, v) \,=\, h^E(\nabla u, v) + h^E(u, \nabla v), \qquad u, v
\in C^{\infty}(M, E)
\]
and $h^E$ can be extended canonically to a sesquilinear map
\[
h^E \, : \, \Omega^{\bullet}(M, E) \times \Omega^{\bullet}(M, E) \rightarrow \Omega^{\bullet}(M, \mathbb{C}).
\]
For $\omega_1, \omega_2 \in \Omega^{\bullet}(M, E)$ and for each $j = 0,\cdots,d,$ we then obtain a sesquilinear pairing
\begin{equation}\label{E:paring}
h^E \, : \, \Omega^{j}(M, E) \times \Omega^{d-j}(M, E) \rightarrow \mathbb{C}, \quad (\omega_1, \omega_2) \mapsto \int_M h^E(\omega_1 \wedge \omega_2).
\end{equation}
The pairing (\ref{E:paring}) induces a non-degenerate sesquilinear pairing
\begin{equation}
H^j(\,M, E\,) \otimes H^{d-j}(\,M, E\,) \longrightarrow \mathbb{C}, \qquad j=0\cdots d,
\end{equation}
and allows us to identify $H^j(\,M, E\,)$ with the dual space of $H^{d-j}(\,M, E\,)$.
Using the construction of Subsection~3.4 of \cite{BK3} we thus obtain a canonical involution
\begin{equation}\label{E:invol}
    D : \, \operatorname{Det}\big(H^{\bullet}(M,E)\big) \ \longrightarrow \ \operatorname{Det}\big(H^{\bullet}(M,E)\big).
\end{equation}

Note that if the flat bundle $E$ is acyclic, then the complex determinant line $\operatorname{Det}\big(H^{\bullet}(M,E)\big)$ is canonically isomorphic to $\mathbb{C}$ and under this isomorphism the involution \eqref{E:invol} coincides with the complex conjugation.


If $h \in \operatorname{Det}\big(H^{\bullet}(M,E)\big)$ and $D(h) \,=\, h$, then the element $h$ will be called {\em real} . The real elements of $\operatorname{Det}\big(H^{\bullet}(M,E)\big)$
form a real line.

If $h \in \operatorname{Det}\big(H^{\bullet}(M,E)\big)$ can be represented in the form $h \,=\, h_0e^{i\phi}$, where $h_0$ is real, then $\phi \in \mathbb{R}$ will be called the {\em phase} of $h$. It is defined up to an integral multiple of  $\pi$ and we will denote it by $\operatorname{Ph}(h)$.

\subsection{On sign conventions}\label{S:convention}
The definition of the canonical involution $D$ in \cite[Subsection
10.1]{BK3} is different from the definition of the canonical involution in \cite[Subsection 2.1]{FT1} by a factor $(-1)^{\nu},\, \nu \in \mathbb{Z}$. Hence, if we denote by $\widetilde{\operatorname{Ph}}(h)$ the phase of $h \in \operatorname{Det}\big(H^{\bullet}(M,E)\big)$ as it is defined in \cite{FT1}, then
\begin{equation}\label{E:signcon}
\operatorname{Ph}(h) \,=\, \widetilde{\operatorname{Ph}}(h) + \frac{\pi \nu}{2} \qquad \operatorname{mod} \pi\mathbb{Z}.
\end{equation}
Note, however, that if the bundle $E$ is acyclic, then both involutions coincide with the complex conjugation, cf. Subsection 3.1 and ~\cite[Lemma 2.2]{FT1}. Hence, for the acyclic case, $\nu \,=\, 0$ and $\operatorname{Ph}(h)\,=\,\widetilde{\operatorname{Ph}}(h)$.

\subsection{Phase of the refined analytic torsion.}
In this subsection we derive the formula of the phase of the refined
analytic torsion $\rho_{\operatorname{an}}(\nabla)$. We have the
following proposition.
\begin{proposition}\label{P:phaseRAT}
Let $M$ be a closed oriented manifold of odd dimension $d=2n-1$ and
let $E$ be a flat complex vector bundle over $M$ admitting a flat
Hermitian metric and endowed with a flat connection $\nabla$. Then
the phase of the refined torsion $\rho_{\operatorname{an}}(\nabla)$
is given by the following formula:
\begin{equation}
\operatorname{Ph}\big(\,\rho_{\operatorname{an}}(\nabla)\,\big) \, =
\, -\pi\big(\,\eta (\,\nabla) - \operatorname{rank}E \cdot
\eta_{\operatorname{trivial}}\,\big)\, \qquad \operatorname{mod} \pi\mathbb{Z}.
\end{equation}
\end{proposition}
\begin{proof}
From \cite[Theorem 10.3]{BK3}, we have
\begin{equation}
D\big(\,\rho_{\operatorname{an}}(\nabla)\,\big) \,=\,
\rho_{\operatorname{an}}(\nabla) \cdot e^{2i\pi
\big(\,\eta(\nabla) - \operatorname{rank}E \cdot
\eta_{\operatorname{trivial}}\,\big)}.
\end{equation}
If $\rho_{\operatorname{an}}(\nabla) \,=\, \rho_0e^{i\phi}$, where
$\rho_0$ is real with respect to the canonical involution
\eqref{E:invol}, then we obtain
$D\big(\,\rho_{\operatorname{an}}(\nabla)\,\big)/\rho_{\operatorname{an}}(\nabla)
\,=\, e^{-2i\phi}$. Therefore,
\begin{equation}
\operatorname{Ph}\big(\,\rho_{\operatorname{an}}(\nabla)\,\big) \, =
\, -\pi \big(\,\eta(\nabla) - \operatorname{rank}E \cdot
\eta_{\operatorname{trivial}}\,\big)\, \qquad \operatorname{mod} \pi\mathbb{Z}.
\end{equation}
\end{proof}

\subsection{Phase of the Farber-Turaev torsion.}

The homological version of the formula of the phase of the
Farber-Turaev torsion was computed in \cite{FT1}. Similarly we have
cohomological version of the formula of the phase of the
Farber-Turaev combinatorial torsion $\rho_{\varepsilon,o}(\nabla)$.

Following Farber~\cite{F}, we denote by
$\operatorname{Arg}_\nabla$ the unique cohomology class
$\operatorname{Arg}_\nabla \in H^1(M,\mathbb{C}/\mathbb{Z})$ such
that for every closed curve $\gamma \in M$ we have
\begin{equation}\label{E:Mon}
    \det\big(\, \operatorname{Mon}_{\nabla}(\gamma)\, \big) \ = \  \exp\big( \,2\pi i \langle\, \operatorname{Arg}_\nabla,[\gamma] \,\rangle\, \big),
\end{equation}
where $\operatorname{Mon}_\nabla(\gamma)$ denotes the monodromy of
the flat connection $\nabla$ along the curve $\gamma$ and
$\langle\cdot,\cdot\rangle$ denotes the natural pairing
\[
    H^1(M,\mathbb{C}/\mathbb{Z}) \times H_1(M,\mathbb{Z})\  \longrightarrow \ \mathbb{C}/\mathbb{Z}.
\]
Note that when $\nabla$ is a Hermitian connection,
$\operatorname{Mon}_{\nabla}(\gamma)$ is unitary and
$\operatorname{Arg}_\nabla \in H^1(M,\mathbb{R}/\mathbb{Z})$, cf.
\cite{F}.

Let $c(\varepsilon) \in H_1(M, \mathbb{Z})$ denote the
characteristic class of the Euler structure $\varepsilon$, cf. \cite[Subsection 5.2]{FT2} or \cite{T2}, then we have the following
proposition, \cite[Theorem 2.3]{FT1}.

\begin{proposition}\label{P:phaseFT}
Let $M$ be a closed oriented manifold of odd dimension $d=2n-1$ and
let $E$ be a flat complex vector bundle over $M$ admitting a flat
Hermitian metric and endowed with a flat connection $\nabla$. Then
the phase of the Farber-Turaev torsion
$\rho_{\varepsilon,o}(\nabla)$ is given by the following formula:
\begin{equation}
\operatorname{Ph}\big(\,\rho_{\varepsilon,o}(\nabla)\,\big) \, = \, \pi \langle\, \operatorname{Arg}_\nabla, c(\varepsilon) \,\rangle \,+\, \frac{\pi\nu}{2}  \qquad
\operatorname{mod} \pi\mathbb{Z},
\end{equation}
where $\nu \in \mathbb{Z}$. If, moreover, the bundle $E$ is acyclic, then $\nu \,=\, 0$, cf. Subsection~\ref{S:convention}.
\end{proposition}

\subsection{Comparison between the Farber-Turaev and the refined analytic torsions.}\label{T:trat}
In \cite{BK1,BK4}, Braverman and Kappeler computed the ratio
\[
    R \ = \ R(\nabla,\varepsilon,o) \ := \ \frac{\rho_{\operatorname{an}}(\nabla)}{\rho_{\varepsilon,o}(\nabla)}.
\]
We now briefly remind their result. First, we need to introduce some
additional notations.

Let us denote by $\widehat{L}(p) \in H_{\bullet}(M, \mathbb{Z})$ the
Poincar\'{e} dual of the cohomology class $[L(p)]$, where $L(p) =
L_M(p)$ is the Hirzebruch $L$-polynomial in the Pontrjagin forms of
the Riemannian metric $g^M$. Let $\widehat{L}_1 \in
H_1(M,\mathbb{Z})$ denote the component of $\widehat{L}(p)$ in
$H_1(M,\mathbb{Z})$. Then
\[
\langle \,[L(p)] \cup \operatorname{Arg}_{\nabla},
[M]\,\rangle \  = \  \langle\, \operatorname{Arg}_{\nabla},
\widehat{L}_1\,\rangle \ \in\ \mathbb{C}/\mathbb{Z}.
\]
Note that  when $\dim{M} \equiv 3(\Mod 4), \widehat{L}_1 =0.$

The following Braverman-Kappeler theorem comparing the refined
analytic torsion with Farber-Turaev torsion was proven in
\cite[Theorem 5.11]{BK4}. We will restrict to the case that the
connected component $C$ of $\operatorname{Flat}(E)$ contains a
Hermitian connection.

\begin{theorem}\label{T:ratio}
Suppose that $M$ is a closed oriented odd dimensional manifold. Let
$E$ be a flat complex vector bundle over $M$ admitting a Hermitian
metric and endowed with a flat connection $\nabla$. Let
$\varepsilon$ be an Euler structure on $M$ and let $o$ be a
cohomological orientation of $M$. Then, for each connected component
$C$ of the set $\operatorname{Flat}(E)$ that contains a Hermitian
connection, there exists a constant $\theta^C = \theta^C_o \in
\mathbb{R}/2\pi\mathbb{Z}$, depending on $o$ (but not on
$\varepsilon$), such that, for any connection $\nabla \in C$,
\begin{equation}\label{E:ratio}
\frac{\,\rho_{\operatorname{an}}(\nabla)}{\rho_{\varepsilon,o}(\nabla)}
\,=\, \pm e^{i\theta^C_o} \cdot e^{-\pi
i\langle\,\operatorname{Arg}_{\nabla} ,
c(\varepsilon)+\widehat{L}_1\,\rangle}.
\end{equation}
\end{theorem}

Now we compute the constant $\theta^C$ which appears in the quotient
of the refined analytic torsion and the cohomological Farber-Turaev
torsion of $M$, cf. \eqref{E:ratio}. We have the following theorem.

\begin{theorem}\label{T:const}
Suppose that $M$ is a closed oriented odd dimensional manifold. Let
$E$ be a flat complex vector bundle over $M$ admitting a Hermitian
metric and endowed with a flat connection $\nabla$. Let
$\varepsilon$ be an Euler structure on $M$ and let $o$ be a
cohomological orientation of $M$. If the connected component $C$ of
$\operatorname{Flat}(E)$ contains a Hermitian connection, then, for some $\nu \in \mathbb{Z}$,
\begin{enumerate}
\item
If $\dim{M} \equiv 1(\Mod 4)$, then
\begin{equation}\label{E:constant1}
\theta^C \ = \ - \pi \big(\, \operatorname{Re}\eta(\nabla)  -
\operatorname{Re} \langle\, \operatorname{Arg}_{\nabla},
\widehat{L}_1 \,\rangle \,\big)\,+\, \frac{\pi\nu}{2}  \qquad \Mod \quad \pi\mathbb{Z}.
\end{equation}
\item
If $\dim{M} \equiv 3(\Mod 4)$, then
\begin{equation}\label{E:constant2}
\theta^C \ = \  - \pi \big(\, \operatorname{Re}\eta(\nabla) -
\operatorname{rank}E \cdot \eta_{\operatorname{trivial}}\,\big)\,+\, \frac{\pi\nu}{2}
\qquad \Mod \quad \pi\mathbb{Z}.
\end{equation}
\end{enumerate}
If, moreover, the bundle $E$ is acyclic, then $\nu \,=\, 0$, cf. Subsection~\ref{S:convention}.
\end{theorem}

\begin{proof}
If $\nabla \in C$ is a Hermitian connection, then the theorem
follows by combining Proposition~\ref{P:phaseRAT} and
Proposition~\ref{P:phaseFT} with Theorem~\ref{T:ratio}.

Suppose that $\nabla_t (t \in [0,\/1])$ is a smooth family of
connections in $C$ such that $\nabla_0\,=\,\nabla$ is Hermitian. From Theorem 12.3 and Lemma 12.6 of \cite{BK1} we
conclude that
\begin{equation}\label{E:eq1}
\frac{d}{dt}\eta(\nabla_t) \,=\, \frac{d}{dt}\langle
\operatorname{Arg}_{\nabla_t}, \widehat{L}_1 \rangle .
\end{equation}
Lemma 5.5 of \cite{BK4} shows that
\begin{equation}\label{E:eq2}
\exp \big(\, \pi \operatorname{Im} \langle\, \operatorname{Arg}_{\nabla},
\widehat{L}_1 \,\rangle \,\big) \,=\, \exp \big(\, \pi \operatorname{Im} \eta(\nabla) \,\big).
\end{equation}
Hence by combining \eqref{E:eq1}, \eqref{E:eq2} with the case that
$\nabla \in C$ is a Hermitian connection, the theorem follows.
\end{proof}


From Theorem~\ref{T:ratio} and Theorem~\ref{T:const}, we have the
following theorem which improves the Braverman-Kappeler theorem
(Theorem~\ref{T:ratio}).

\begin{theorem}\label{T:comp}
Suppose that $M$ is a closed oriented odd dimensional manifold. Let
$E$ be a flat complex vector bundle over $M$ admitting a Hermitian
metric and endowed with a flat connection $\nabla$. Let
$\varepsilon$ be an Euler structure on $M$ and let $o$ be a
cohomological orientation of $M$. If the connected component $C$ of
$\operatorname{Flat}(E)$ contains a Hermitian connection, then
\begin{equation}\label{E:ratt}
\frac{\,\rho_{\operatorname{an}\,}(\nabla)}{\rho_{\varepsilon,o}(\nabla)\,}
\,=\, \pm i^{\nu} \cdot
e^{-\pi i\langle\,\operatorname{Arg}_{\nabla},c(\varepsilon) \, \rangle}
\cdot e^{-i\pi(\eta(\nabla)- \operatorname{rank}E
\cdot \eta_{\operatorname{trivial}})},
\end{equation}
where $\nu \in \mathbb{Z}$.
If, moreover, the bundle $E$ is acyclic, then $\nu \,=\, 0$, cf. Subsection~\ref{S:convention}.\end{theorem}


\section{Comparison between the Farber-Turaev and the Burghelea-Haller torsions.}

In \cite{BH1, BH2} Burghelea and Haller introduced a refinement of
the square of the Ray-Singer torsion for a closed manifold of
arbitrary dimension, provided that the complex vector bundle $E$ admits
a non-degenerate complex valued symmetric bilinear form $b$. They
defined a complex valued quadratic form $\tau_{b, \nabla}$ on the
determinant line $\operatorname{Det}\big(H^{\bullet}(M,E)\big)$.
Then they defined a complex valued quadratic form, referred to as
{\em complex Ray-Singer torsion}. For the closed oriented odd
dimensional manifold $M$ and the complex vector bundle $E$ over $M$
endowed with a flat connection $\nabla$, it is given by
\begin{equation}
\tau^{\operatorname{BH}}_{b, \gamma, \nabla} \,:=\, \tau_{b, \nabla} \cdot
e^{-2\int_M \omega_{\nabla, b}\wedge \gamma},
\end{equation}
where $\gamma \in \Omega^{d-1}(M)$ is an arbitrary closed
$(d-1)$-form and $\omega_{\nabla, b} \in \Omega^1(M)$ is the
Kamber-Tondeur form, cf. \cite[Section 2]{BH2}.

Burghelea and Haller conjectured, \cite[Conjecture 5.1]{BH2}, that
for a suitable choice of $\gamma$ the form $\tau^{\operatorname{BH}}_{b, \gamma,
\nabla}$ is roughly speaking equal to the square of the
Farber-Turaev torsion and established the conjecture in some
non-trivial situations. Though the conjecture is for manifolds of
arbitrary dimensions, we restrict to the odd dimensional case and
adopt the following formulation from \cite[Conjecture 1.9]{BK5}.

\begin{conjecture}{\rm\textbf{[Burghelea-Haller]}}\label{C:conj}
Let $M$ be a closed oriented manifold of odd dimension $d=2n-1$ and
let $E$ be a flat complex vector bundle over $M$ endowed with a flat
connection $\nabla$. Let $b$ be a non-degenerate symmetric bilinear
form on $E$. Let $\varepsilon$ be an Euler structure on $M$
represented by a non-vanishing vector field $X$ and let $o$ be a
cohomological orientation of $M$. Fix a Riemannian metric $g^M$ on
$M$ and let $\Psi(g^M) \in \Omega^{d-1}(TM\backslash\{0\})$ denote
the Mathai-Quillen form, \cite[pp. 40-44]{BismutZhang},
\cite[section7]{MathaiQuillen}. Set
\begin{equation}
\gamma_{\varepsilon} \,=\, \gamma_{\varepsilon}(g^M) \,:=\,
X^{\ast}\Psi(g^M).
\end{equation}
 Then
\begin{equation}
\tau^{\operatorname{BH}}_{b, \gamma_{\varepsilon}
,\nabla}(\rho_{\varepsilon,o}(\nabla)) \,=\, 1.
\end{equation}
\end{conjecture}

In \cite{BK5} Braverman and Kappeler expressed the Burghelea-Haller
complex Ray-Singer torsion in terms of the square of the refined
analytic torsion $\rho_{\operatorname{an}}(\nabla)$ and the eta
invariant $\eta(\nabla)$. In particular, they proved the following weak
version of the Burghelea-Haller conjecture: $\tau^{\operatorname{BH}}_{b, \gamma_{\varepsilon}
,\nabla}(\rho_{\varepsilon,o}(\nabla))$ is a locally constant in $\nabla$ and
\[
|\tau^{\operatorname{BH}}_{b, \gamma_{\varepsilon}
,\nabla}(\rho_{\varepsilon,o}(\nabla))| \,=\, 1.
\]
In the following theorem we improve this result for the case when
$\nabla$ belongs to a connected component of the space of flat
connections on the associated complex vector bundle $E$ which
contains a Hermitian connection. More precisely, we have the
following theorem:

\begin{theorem}\label{T:thm}
Under the assumptions of Conjecture~\ref{C:conj} and assume that the
connected component $C$ of the set $\operatorname{Flat}(E)$ of flat
connections on $E$ contains a Hermitian connection, then
\begin{equation}\label{E:equality}
\tau^{\operatorname{BH}}_{b, \gamma_{\varepsilon}
,\nabla}(\rho_{\varepsilon,o}(\nabla)) \,=\, \pm 1, \qquad  for \,
all \quad \nabla \in C.
\end{equation}
\end{theorem}

\begin{proof}
In Theorem 1.10 of \cite{BK5}, Braverman and Kappeler proved that
$\tau^{\operatorname{BH}}_{b, \gamma_{\varepsilon}
,\nabla}(\rho_{\varepsilon,o}(\nabla))$ is constant on $C$ with
absolute value 1. Hence it is enough to prove the equality
\eqref{E:equality} in the case when $\nabla \in C$ is a Hermitian
connection.

From (1.1), (1.2) and Theorem 1.4 of \cite{BK5}, we have
\begin{equation}\label{E:bkbh}
\tau^{\operatorname{BH}}_{b, \gamma_{\varepsilon}
,\nabla}(\rho_{\operatorname{an}}(\nabla))\,=\,\pm e^{-2\pi
i\big(\,\eta(\nabla)-\operatorname{rank}E\cdot\eta_{\operatorname{trivial}}
\,\big)}\cdot e^{-2\int_M\omega_{\nabla,
b}\wedge\gamma_{\varepsilon}}.
\end{equation}
We also have, cf. \cite[Subsection 5.4]{BK5},
\begin{equation}\label{E:mono}
e^{-2\int_M\omega_{\nabla, b}\wedge\gamma_{\varepsilon}} \,=\,
e^{<[\omega_{\nabla,b},c(\varepsilon)]>}\,=\,\pm
e^{-2\pi i\langle\,\operatorname{Arg}_{\nabla},c(\varepsilon) \, \rangle}.
\end{equation}
From \eqref{E:Mon} and \eqref{E:ratt}, we get
\begin{equation}\label{E:ratt1}
\Big(\,\frac{\,\rho_{\operatorname{an}\,}(\nabla)}{\rho_{\varepsilon,o}(\nabla)\,}\,\Big)^2
\,=\, \pm
e^{-2\pi i\langle\,\operatorname{Arg}_{\nabla},c(\varepsilon) \, \rangle}
\cdot e^{-2i\pi\big(\,\eta(\nabla)- \operatorname{rank}E \cdot
\eta_{\operatorname{trivial}}\,\big)}.
\end{equation}
By combining \eqref{E:bkbh} and \eqref{E:mono} with \eqref{E:ratt1}, we
obtain the result. \hfill$\square$
\end{proof}


\section{Refined analytic torsion of lens spaces}
In this section we compute the refined analytic torsion of lens
spaces. We begin with recalling the relationship of the acyclic case
of the refined analytic torsion with the Ray-Singer torsion and the
eta invariants. We then recall the definition of a lens space and the
formula for the Ray-Singer torsion of a lens space from \cite{R}.
Then we recall the formula for the eta invariant of a lens space
from \cite{APS2}. By combining these results, we obtain the refined
analytic torsion of a lens space.

\subsection{The acyclic case of refined analytic torsion}\label{SS:twisted}
Denote by $\widetilde{M}$ the universal covering of $M$ and by
$\pi_1(M)$ the fundamental group of $M$, viewed as the group of deck
transformations of $\widetilde{M} \rightarrow M.$ For each complex
representation $\alpha : \pi_1(M) \rightarrow GL(r, \mathbb{C})$, we
denote by
\begin{equation}\label{E:Ealp}
E_{\alpha} \ := \  \widetilde{M} \times_{\alpha} \mathbb{C}^r
\rightarrow M
\end{equation}
the flat vector bundle induced by $\alpha$. Let $\nabla_{\alpha}$ be
the flat connection on $E_{\alpha}$ induced from the trivial
connection on $\widetilde{M} \times \mathbb{C}^r$. We also denote by
$\nabla_{\alpha}$ the induced differential
\[
\nabla_{\alpha}  \ : \  \Omega^{\bullet}(\,M, \,E_{\alpha}\,)
\rightarrow  \Omega^{\bullet + 1}(\,M,\, E_{\alpha}\,),
\]
where $\Omega^{\bullet}(M, E_{\alpha})$ denotes the space of smooth
differential forms of $M$ with values in $E_{\alpha}$.

If the representation $\alpha$ is acyclic, i.e.,
$H^{\bullet}(M,E_{\alpha}) = 0$, then the determinant line
$\text{Det}(H^{\bullet}(M,E_{\alpha}))$ is canonically isomorphic to
$\mathbb{C}$. In particular, if $\alpha$ is an acyclic unitary
representation of $\pi_1(M)$, then, cf.~\cite[Section 12]{BK1},
\begin{equation}\label{E:relation}
\rho_{\text{an}}(\nabla_{\alpha}) \ = \ \rho^{\operatorname{RS}}_{\alpha} \cdot
e^{-i\pi\eta_{\alpha}}\cdot e^{i\pi \cdot \rk \alpha \cdot
\eta_{\operatorname{trivial}}},
\end{equation}
where $\rho^{\operatorname{RS}}_{\alpha} := \rho^{\operatorname{RS}}(\nabla_{\alpha})$ is the
well-known Ray-Singer torsion, \cite{RS}, and $\eta_{\alpha} :=
\eta(\nabla_{\alpha})$.

\subsection{The lens space}
Fix an integer $m\ge3$ and let $G_m$ denote the cyclic group of order $m$. We fix a generator $g\in G_m$ so that $G_m= \{1,g,g^2,\ldots,g^{m-1}\}$.

Let $p_1,\ldots,p_n$ be integers relatively prime to $m$. Then the action of $G_m$ on the sphere
\[
    S^{2n-1}\ =\ \big\{\,z\in \mathbb{C}^n\mid\,\|z\|=1\,\big\}
\]
defined by
\begin{equation}\label{E:action}
    g\cdot (\,z_1,\ldots,z_n\,) \ = \ \big(\,e^{2\pi i p_1/m}z_1,\ldots,e^{2\pi i p_n/m}z_n\,\big)
\end{equation}
is free. The {\em lens space} $L=L(m;p_1,\ldots,p_n)$ is the orbit space of this action
\[
    L\ = \ L(\,m;p_1,\ldots,p_n\,)\ := \ S^{2n-1}/G_m.
\]
Clearly, $\pi_1(L)= G_m$.

Fix $q\in \mathbb{Z}$ and consider the unitary representation $\alpha=\alpha_q:\pi_1(L)= G_m \to U(1)$, defined by
\[
    \alpha_q(g) \ = \ e^{2\pi i q/m}.
\]
We will be interested in the refined analytic torsion
$\rho_{\text{an}}(q)=\rho_{\text{an}}(\nabla_{\alpha_q})$ associated
to the representation $\alpha_q$.


\subsection{The Ray-Singer torsion of the lens spaces}

In this subsection we recall the formula for the Ray-Singer torsion of lens spaces from ~\cite[Section 4]{R}.
Note that our definition of logarithm of Ray-Singer torsion is
negative one half of the logarithm of the Ray-Singer torsion in \cite{R}.

\begin{proposition}
Let $l_k$$(k=1...n)$ be any integers such that $l_kp_k\equiv1 (\Mod
m)$ and let $\rho^{\operatorname{RS}}_{\alpha_q}(L)$ denote the Ray-Singer torsion
of the lens space $L$ associated to the nontrivial acyclic
representation $\alpha_q$, then
\[
\rho^{\operatorname{RS}}_{\alpha_q}(L) \ =\ \prod_{k=1}^n\,  \mid e^{\frac{2\pi i q
l_k}{m}}-1\mid.
\]

\end{proposition}

\subsection{The $\eta$ invariant of the odd signature operator for lens spaces}

A slight modification of Proposition 2.12 in ~\cite{APS2}, where the eta invariant $\eta_{\operatorname{trivial}}$ of a lens space for trivial representation was
computed, we have the following proposition, see also~\cite[Proposition 4.1]{D}.

\begin{proposition}\label{P:etaq}
Let $\eta_{\alpha_q}  =  \eta_q$ denote the eta invariant of the odd signature operator $\B_{\even}(\nabla_{\alpha_q}, g^L)$ of the lens space $L$, then
\[
  \eta_q \ = \ \frac{i^{-n}}{2m}\,\sum_{l=1}^{m-1}\,\Big(\,e^{\frac{2\pi i l q}{m}}\cdot\prod_{j=1}^n\cot\frac{\pi l p_j}{m} \,\Big).
\]
In particular, when $n$ is even we have
\[
 \eta_q \ = \ \frac{i^{-n}}{2m}\sum_{l=1}^{m-1}\,\Big(\,\cos\frac{2\pi l q}{m}\cdot\prod_{j=1}^n\cot\frac{\pi l p_j}{m}\,\Big)
\]
and when $n$ is odd we have
\[
 \eta_q \ = \ \frac{i^{1-n}}{2m}\sum_{l=1}^{m-1}\Big(\,\sin\frac{2\pi l q}{m}\cdot\prod^n_{j=1}\cot\frac{\pi l p_j}{m}\,\Big).
\]
\end{proposition}

Note that our $\eta$ invariant is equal to one half of the $\eta$ invariant in \cite{APS2}.

Combined these two propositions with ~\eqref{E:relation} and Proposition 2.12 in ~\cite{APS2}, we have
\begin{theorem}\label{T:analtorsion}
For the lens space $L\ = \ L(\,m;p_1,\ldots,p_n\,), m \ge 3$ and the nontrivial acyclic representation $\alpha_q:\pi_1(L) \ = \ G_m \to U(1)$ such that $\alpha_q(g) \ = \ e^{2\pi i q/m}$, the refined analytic torsion
\[
\rho_{\operatorname{an}}(q) \ = \ \prod_{k=1}^n \mid e^{\frac{2\pi i
q l_k}{m}}-1\mid \cdot e^{
\frac{i^{-(n+1)}\pi}{2m}\sum_{l=1}^{m-1}(e^{\frac{2\pi i l
q}{m}}\cdot\prod_{j=1}^n\cot\frac{\pi l p_j}{m})} \cdot e^{
\frac{i^{-(n-1)}\pi}{2m}\sum_{l=1}^{m-1}(\prod_{j=1}^n\cot\frac{\pi
l p_j}{m})},
\]
where $l_k$$(k=1,\cdots,n)$ are any integers such that $l_kp_k\equiv1 (\Mod m)$

\end{theorem}


\section{Comparison of the refined analytic and the Turaev torsions of a lens space}

In this section we begin with recalling the definition of the Turaev torsion and computing the Turaev torsion of lens spaces.  We then calculate the constant $\theta^C$ for the three-dimensional lens space $L(5; 1,1)$ and
the ratio $R$ of the refined analytic torsion and cohomological Turaev torsion of the five-dimensional lens space $L(3; 1,1,1)$ and explain how our computations give
answers to Questions~1 and 2 of the introduction.

\subsection{Torsion of an acyclic chain complex}\label{SS:cc}

Let $\mathbb{F}$ be a field of characteristic zero and let
\[
    C:\quad 0 \ {\longrightarrow}\ C_d\ \stackrel{\partial_{d-1}}{\longrightarrow}\ C_{d-1}\ \stackrel{\partial_{d-2}}{\longrightarrow}\ \cdots\
    \stackrel{\partial_0}{\longrightarrow}\ C_0\ {\longrightarrow}\ 0
\]
be a finite dimensional chain complex over $\mathbb{F}$. Assume
that the chain complex $(C,\partial)$ is acyclic, i.e.
$H_{\ast}(C)=0$. For each $i$, let $c_i$ be a fixed basis for $C_i$
and  $b_i$ be a sequence of vectors in $C_i$ whose image under
$\partial_{i-1}$ is a basis in $\operatorname{Im} \partial_{i-1}$.
Then the vectors $\partial_i(b_{i+1}), b_i$ form a basis for $C_i$.
The torsion of  $C$ is defined by
\[
\tau_{\text{comb}}(C) \ = \  \prod_{i=0}^d\, [\,
\partial_i(b_{i+1})b_i/c_i \,]^{(-1)^{i+1}},
\]
where $[\partial_i(b_{i+1})b_i/c_i]$ is the determinant of the
matrix transforming $c_i$ into the basis $\partial_i(b_{i+1}), b_i$
of $C_i$.


\subsection{The Reidemeister torsion}\label{SS:reid}

Fix a CW-decomposition $X\,=\,\{e_1,\ldots,e_N\}$ of $M$. For each
$j=1,\ldots,N$, fix a lift $\widetilde{e}_j$, i.e. a cell of the
CW-decomposition $\widetilde{X}$ of $\widetilde{M}$, such that
$\pi(\widetilde{e}_j)= e_j$. By \eqref{E:Ealp}, the pull-back of the
bundle $E_{\alpha}$ to $\widetilde{M}$ is the trivial bundle
$\widetilde{M}\times \ \mathbb{C}^n\to \widetilde{M}$. Hence, the set of
the cells $\widetilde{e}_1\ldots \widetilde{e}_N$ identifies the chain
complex $C(X,\alpha)$ of the CW-complex $X$ with coefficients in
$E_{\alpha}$ with the complex $C(\widetilde{X}) \otimes_{\alpha} \mathbb{C}^r$, where $\alpha : \pi_1(X) \rightarrow GL(r, \mathbb{C})$ is a representation. Assume that this chain complex $C(X,\alpha)$ is acyclic, i.e.
\[
H_{\ast}(\,M,E_{\alpha}\,) \ = \ H_{\ast}(\,C(X,\alpha)\,) \ = \
0,
\]
then the Reidemeister torsion is defined as the torsion of this
chain complex.


\subsection{Combinatorial Euler structures and homological Turaev torsion }\label{SS:Euler}
In this subsection we recall the definition of combinatorial Euler
structures from ~\cite{T3}.

A family $\hat{e}=\{\hat{e}_i\}$ of open cells in the maximal
abelian covering
\[
\widehat{X} \ = \ \widetilde{X}/[\,\pi_1(X), \pi_1(X)\,]
\]
of $X$ is called {\em fundamental} if each open cell $e_i$ in $X$ is
covered exactly by one cell $\hat{e}_i$ of $\hat{e}$.

Following Turaev, we denote the operation of any two cells in
multiplicative notation. Let
\[
\hat{e}'/\hat{e} \ = \ \prod_{e_i \in
X}\,(\,\hat{e}'_i/\hat{e}_i\,)^{(-1)^{\dim e_i}} \in H_1(M)
\]
for any two fundamental families $\hat{e}$ and $\hat{e}'$, here
$\hat{e}'_i/\hat{e}_i \in H_1(M)$.  We say that the fundamental
families $\hat{e}$ and $\hat{e}'$ are equivalent if $
\hat{e}/\hat{e}' \,=\, 1$. The equivalence classes are called {\em
combinatorial Euler structures} on $M$.

Let $\alpha : \pi_1(M) \rightarrow GL(r, \mathbb{C})$ be an acyclic
representation. Then we can associate each combinatorial Euler
structure $\varepsilon$ on $M$ the homological Turaev torsion
\[
\tau_{\alpha}(\,M,\varepsilon\,) \,=\, \tau_{\alpha}(\,M,\hat{e}\,)
\in \mathbb{C} / \pm.
\]


For each Euler structure $\varepsilon$ on $M$, there is an {\em
Euler class} $c(\varepsilon)\in H_1(M)$ associated to it, cf.
\cite{T2} or \cite[Subsection 5.2]{FT2}. If $d = \dim M$ is odd,
then (in multiplicative notation)
\begin{equation}
c(h\varepsilon) = h^2c(\varepsilon)
\end{equation}
for any $\varepsilon \in \Eul(M)$, $h \in H_1(M)$.

Turaev also introduced the homology orientation to get rid of the
sign indeterminacy of Reidemeister torsion. For our purpose it will
be enough to consider the Turaev torsion up to sign, so we skip the
definition of homology orientation.


\subsection{Cohomological Turaev torsion}\label{SS:cohtor}
In this subsection we recall the relationship of the cohomological
Turaev torsion with the homological Turaev torsion, cf.
\cite[Subsection 9.2]{FT2}.

Let $M$ be a closed oriented manifold of odd dimension $d=2n-1$,
where $n \ge 1$ and $\alpha : \pi_1(M) \rightarrow GL(r,
\mathbb{C})$ be an acyclic representation of the fundamental group
of $M$. Denote the cohomological Turaev torsion associated to the Euler structure
$\varepsilon$ by $\rho_{\alpha}(\,M,\varepsilon\,)(\, = \rho_{\varepsilon,o}(\nabla_{\alpha})\,)$, then, cf. \cite[p. 218 (9.2), (9.3)]{FT2},
\begin{equation}\label{E:cohtorsion}
\rho_{\alpha}(\,M,\varepsilon\,) \ = \
\frac{1}{\tau_{\alpha^{\ast}}(\,M,\varepsilon\,)},
\end{equation}
where $\alpha^{\ast}$ is the dual representation of $\alpha$. Recall
that, for all $g \in \pi_1(M)$, $\alpha^{\ast}(g) \,=\,
(\,\alpha(g)^{-1}\,)^t$, cf. \cite[subsection 4.1]{FT2}, where $t$
denotes the transpose of matrices. It is clear that for all $g \in
\pi_1(M)$ we have $\det(\,\alpha(g)\,) \cdot
\det(\,\alpha^{\ast}(g)\,) \,=\, 1$.

\subsection{The Turaev torsion of lens spaces }
In this subsection we compute the Turaev torsion of lens spaces. Let
$L=L(m; p_1,\cdots,p_n),m \ge 3$, be the lens space. First we fix a preferred Euler structure $\epsilon$ on $L$. Consider the
CW-decomposition  $e =  \{\,e_j\,\}_{j \ = \ 1, \ldots, 2n-1}$ of $L$ such that the
CW-decomposition $e$ lifts to a $G_m$-equivariant CW-decomposition of
$S^{2n-1}$. More precisely, for each $j = 1, \ldots, 2n-1$, let us fix the
lift $\widetilde{e}_j$ of $e_j$ to $S^{2n-1}$ such that, for each $i = 1, \ldots, n$,
\begin{equation}\label{E:cell1}
   \widetilde{e}_{2i-1} \ = \ \big\{\,(z_1,\cdots,z_n)\in S^{2n-1} \mid\, z_{i+1} \ = \ \cdots \ = \  z_n \ = \ 0, 0< \arg z_i < 2\pi/m\,\big\}
\end{equation}
and
\begin{equation}\label{E:cell2}
\widetilde{e}_{2i-2}\ = \ \big\{\,(z_1,\cdots,z_n)\in S^{2n-1} \mid\, z_{i+1} \ = \ \cdots \ = \ z_n \ = \ 0, \arg z_i \ = \ 0 \,\big\}.
\end{equation}
Then
\[
  \widetilde{e} \ = \ \big\{\, g^j \cdot \widetilde{e}_{2i-1}, g^j \cdot \widetilde{e}_{2i-2}\, \big\}_{i \ = \ 1 \ldots n, j \in \mathbb{Z}/m\mathbb{Z}}
\]
defines a $G_m$-equivariant CW-decomposition of $S^{2n-1}$ with $m$ cells in each dimension. Note that $e$ has exactly one cell in each dimension. Then by~\eqref{E:action}, we have
\begin{equation}\label{E:action1}
g \cdot \widetilde{e}_{2i-1} \ = \ \big\{\,(z_1,\cdots,z_n)\in S^{2n-1} \mid\, z_{i+1} \ = \ \cdots=z_n \ = \ 0, 2\pi p_i/m< \arg z_i < 2\pi (p_i+1)/m\,\big\}
\end{equation}
and
\begin{equation}\label{E:action2}
g \cdot \widetilde{e}_{2i-2} \ = \ \big\{\,(z_1,\cdots,z_n)\in S^{2n-1} \mid\, z_{i+1} \ = \ \cdots \ = \ z_n \ = \ 0, \arg z_i \ = \  2\pi p_i/m\,\big\}.
\end{equation}
Recall that $S^{2n-1}$ is the universal covering and also the
maximal abelian covering of the lens space $L$, so we can consider
the collection of cells $\hat{e}= \{\widetilde{e}_j\}_{1\le j\le 2n-1}$
in $S^{2n-1}$ as a fundamental family in $S^{2n-1}$. The equivalence
class of this family defines an Euler structure denoted by
$\epsilon$.

In the following proposition we will give the computation of the
homological Turaev torsion $\tau_{\alpha}(L, \epsilon)$ of the lens
space $L$ and the preferred Euler structure $\epsilon$ by using the
same computation of the Reidemeister torsion of lens spaces, see
\cite[Theorem 10.6, p. 45]{T3}
 for the detailed computation of the
Reidemeister torsion of lens spaces.

\begin{proposition}\label{P:turaev}
Let $L=L(m; p_1,\cdots,p_n),m \ge 3$, be the lens space. Let $g \in
\pi_1(L)$ be the generator and let\/ $l_1,\cdots,l_n \in
\mathbb{Z}/m\mathbb{Z}$ such that\/ $l_kp_k \equiv 1(\Mod  m)$.
Let\/ $\alpha_q:\pi_1(L)=G_m\rightarrow U(1)$ be the nontrivial
unitary representation such that $\alpha_q(g)=e^{2\pi i q/m}$, $1
\le q \le m-1$. Also let $\epsilon$ be the Euler structure defined
as above. Then $H_{\ast}(C(L, \alpha_q))=0$ and
\begin{equation}\label{E:Turaev1}
    \tau_{\alpha_q}(\,L,\epsilon\,)\ =\
    \prod_{k=1}^n\,\mid e^{2\pi i q l_k/m}-1 \mid^{-1} \cdot i^n \cdot e^{-\frac{\pi i q \sum_{k=1}^n l_k}{m}} \in \mathbb{C}^{\ast}/\pm.
\end{equation}

\end{proposition}

\begin{proof}
Assume that $\widetilde{e_{i}}$ is oriented for each $i$ such that the
boundary homomorphism of the chain complex $C(S^{2n-1})$ is given
by, cf.~\eqref{E:cell1}, \eqref{E:action2} and recall that $l_ip_i
\equiv 1(\Mod  m)$,
\[
\partial \widetilde{e}_{2i-1} \ = \ g^{l_i}\widetilde{e}_{2i-2}-\widetilde{e}_{2i-2} \  = \  (g^{l_i}-1)\widetilde{e}_{2i-2}
\]
and, cf.~\eqref{E:cell2} and~\eqref{E:action1},
\[
\partial \widetilde{e}_{2i-2} \ = \  \widetilde{e}_{2i-3}+g\widetilde{e}_{2i-3}+\cdots+g^{m-1}\widetilde{e}_{2i-3} \ = \ \sum_{j=0}^{m-1}g^j\widetilde{e}_{2i-3}.
\]
Since, by assumption, $\alpha_q(g)  \ \not= \ 1$ and since $\alpha_q(g)^m \ = \ 1$, we have that
 \[
 \sum_{j=0}^{m-1}\alpha_q(g)^j \ = \  \frac{\alpha_q(g)^m-1}{\alpha_q(g)-1} \ = \ 0.
  \]
Hence we have the chain complex
\[
C(L, \alpha_q)\,=\, C(S^{2n-1}) \otimes_{\alpha_q}\mathbb{C}
=(\cdots \stackrel{0}{\rightarrow} \mathbb{C}\widetilde{e}_{2i-1}
\stackrel{e^{2\pi i q l_i/m}-1}{\longrightarrow}
\mathbb{C}\widetilde{e}_{2i-2} \stackrel{0}{\rightarrow} \cdots).
\]
It is not difficult to see that
\[
H_{\ast}\big(\,C(L, \alpha_q)\,\big)\,=\,0.
\]
It follows from the definitions that
\begin{equation}\label{E:tors}
\tau_{\alpha_q}(\,L,\epsilon\,) \ = \ \pm \prod_{k=1}^n\,(\,e^{2\pi
i q l_k/m}-1\,)^{-1}.
\end{equation}
We now compute the phase $\theta$ of $\tau_{\alpha_q}(L,\epsilon)$. Set $r \ = \ \mid \tau_{\alpha_q}(L,\epsilon) \mid$, then $\tau_{\alpha_q}(L,\epsilon) \ = \ re^{i\theta}$ and
\begin{equation}\label{E:ph1}
{ \tau_{\alpha_q}(\,L,\epsilon\,)}\Big/\,{
{\overline{\tau_{\alpha_q}(L, \epsilon)}}}\ = \  e^{2i\theta}.
\end{equation}
A simple calculation using~\eqref{E:tors} shows that
\begin{equation}\label{E:ph2}
    { \tau_{\alpha_q}(L, \epsilon)}\Big/\,{ {\overline{\tau_{\alpha_q}(L, \epsilon)}}}\ =\
    \prod_{k=1}^n\, \frac{ (\,e^{2\pi i q l_k/m}-1)^{-1}}{ (e^{-2\pi i q l_k/m}-1\,)^{-1}} \ =\ (-1)^n \cdot e^{-2\pi i q \sum_{k=1}^n l_k/m}.
\end{equation}
From \eqref{E:ph1} and \eqref{E:ph2}, we obtain
\[
\theta \ = \ n\pi/2 - \pi q \sum_{k=1}^n l_k/m.
\]
Hence the proposition follows.

\end{proof}

Note that the first component  $\prod_{k=1}^n\mid e^{2\pi i q l_k/m}-1 \mid^{-1}$ of \eqref{E:Turaev1} is the Reidemeister torsion of $L$.


Now we compute the cohomological Turaev torsion of lens spaces. We will follow the same notations as in Proposition \ref{P:turaev}.

\begin{proposition}\label{P:cohturaev}
Let $\rho_{\alpha_q}(L,\epsilon)$ denote the cohomological Turaev
torsion  of the lens space $L$ associated to the preferred Euler
structure $\epsilon$, then
\[
\rho_{\alpha_q}(\,L,\epsilon\,) \ = \ \prod_{k=1}^n\,\mid e^{2\pi i
q l_k/m}-1 \mid \cdot i^n \cdot e^{-\frac{\pi i q \sum_{k=1}^n
l_k}{m}} \in \mathbb{C}^{\ast}/\pm.
\]
\end{proposition}

\begin{proof}

The proposition follows from~\eqref{E:cohtorsion} and
Proposition~\ref{P:turaev}.

\end{proof}


The following theorem gives the cohomological Turaev torsion of a
lens space for an arbitrary Euler structure $\varepsilon$. Recall
that the cardinality of the set of the Euler structures $\Eul(L)$
and the cardinality of $H_1(L)$ are the same and equal to $m$.

\begin{theorem}\label{T:turaev1}
Let $\hat{e}$ be a fundamental family of the preferred Euler
structure $\epsilon$, cf. Proposition~\ref{P:turaev}, and let
$\varepsilon$ be the Euler structure represented by a fundamental
family $\hat{e}'$. Then there exists $s\in\{0, \ldots, m-1 \}$, such
that \/ $\hat{e}'/\hat{e}=g^s$, cf. subsection \ref{SS:Euler}, and
the cohomological Turaev torsion of lens space $L$ associated to the
Euler structure $\varepsilon$ is given by
\begin{equation}\label{E:cohtor1}
\rho_{\alpha_q}(L, \varepsilon)\ =\ \prod_{k=1}^n\,\mid e^{2\pi i q
l_k/m}-1 \mid \cdot i^n \cdot e^{\frac{\pi i q (2s-\sum_{k=1}^n
l_k)}{m}} \in \mathbb{C}^{\ast}/\pm.
\end{equation}

\end{theorem}

\begin{proof}

The theorem follows easily from Proposition~\ref{P:cohturaev} and the following property of the Turaev torsion, cf. ~\cite[(9.4)]{FT2}, that
 \[
\rho_{\alpha_q}(L, \hat{e}') \ = \  \pm
\alpha_q(\,\hat{e}'/\hat{e})\rho_{\alpha_q}(L, \hat{e}\,).
\]
\end{proof}


\subsection{Dependence of the constant $\theta^C$ on the representation. An example}
Theorem~\ref{T:ratio} does not give any information about the
dependence of the constant $\theta^C$ on the connected component
$C$. In this subsection we use the results of previous subsection to
study this dependence in the case of lens spaces. Our goal is to
show that, in general, $\theta^C \ = \ \theta^\alpha$ does depend on
$\alpha$, thus providing a positive answer to Question~1 of the
introduction.

Let $\alpha_q$ be the representation as before and $L=L(5; 1,1)$ be
the lens space. A direct computation using Theorem~\ref{T:const} and
Proposition~\ref{P:etaq} shows that
\[
\theta^{\alpha_1} \ = \ \theta^{\alpha_4} \ = \ -3\pi/10  \qquad
\Mod \quad \pi\mathbb{Z}
\]
and
\[
\theta^{\alpha_2} \ = \ \theta^{\alpha_3} \ = \ -7\pi/10  \qquad
\Mod \quad \pi\mathbb{Z}.
\]
Therefore we conclude that {\em the constant $\theta^{\alpha_q}$
depends on the representation $\alpha_q$}.

\subsection{The ratio of the refined analytic and the Turaev torsions. An example}
It is natural to ask for which representations $\alpha$ one can find
an Euler structure $\varepsilon$ and the cohomological orientation
$o$ such that $\rho_{\text{an}}(\nabla) =
\rho_{\varepsilon,o}(\nabla)$. In this subsection we use
Theorem~\ref{T:analtorsion} and Theorem~\ref{T:turaev1} to show that
the refined analytic torsion and cohomological Turaev torsion of the
five-dimensional lens space $L(3;1,1,1)$ are never equal.

We compute the ratio $R$ of two torsions of the lens space $L(3; 1,1,1)$ for all nontrivial representations (i.e. $q=1,2$) and all Euler structures (i.e. $s = 0, 1, 2$),  see Theorem~\ref{T:turaev1} for the definition of $s$.
A direct computation shows the following table of the ratio $R$.
\smallskip
\begin{center}\setlength{\extrarowheight}{4pt}
\begin{tabular}{|c|c|c|c|}
  \hline
    &$s=0$&$s=1$&$s=2$\\
  \hline
  $q=1$& $R=\pm e^{\frac{5 \pi i}{9}}$ & $R=\pm e^{\frac{8\pi i}{9}}$ &$R=\pm e^{\frac{2 \pi i}{9}}$  \\
  \hline
  $q=2$&  $R=\pm e^{\frac{-5 \pi i}{9}}$ & $R=\pm e^{\frac{ -8 \pi i}{9}}$ &$R=\pm e^{\frac{-2 \pi i}{9}}$ \\
  \hline
\end{tabular}
\end{center}
\medskip

We conclude that {\em for all Euler structures $\varepsilon$ on $L(3; 1,1,1)$ and all representations $\alpha$ of the fundamental group of $L(3; 1,1,1)$, the refined analytic torsion and the Turaev torsion are not equal.} This provides a partial answer to Question~2 of the introduction.


\bibliographystyle{plain}

\end{document}